\documentclass[a4paper,11pt]{article}
\usepackage{amsmath,amscd,amssymb}
\usepackage[frenchb]{babel}
\usepackage{graphicx,psfrag}
\usepackage{xypic}
\input xy
\xyoption{all}

\def \dm {{\textsc{Preuve. }}}

\def \tha #1#2{\noi{\bf#1{\uppercase{\footnotesize{#2}}}}}

\newtheorem{theor}{\tha{T}{h{\'e}or{\`e}me}}[section]
\newenvironment{theo}{
  \begin{theor}\hs -0.2 cm {\bf .} ---  }
{  \end{theor}}

\newtheorem{propo}[theor]{\tha{P}{roposition}}
\newenvironment{prop}{
  \begin{propo}\hs -0.2 cm {\bf .} ---  }
{  \end{propo}}

\newtheorem{lemma}[theor]{\tha{L}{emme}}
\newenvironment{lemm}{
  \begin{lemma}\hs -0.2 cm {\bf .} ---  }
{  \end{lemma}}

\newtheorem{fait}[theor]{\tha{F}{ait}}

\newtheorem{defini}[theor]{\tha{D}{{\'e}finition}}

\newtheorem{corollaire}[theor]{\tha{C}{orollaire}}

\newtheorem{exemple}[theor]{\sc{Exemple}}

\newtheorem{remarq}[theor]{\sc{Remarque}}
\newenvironment{rema}{
  \begin{remarq}\hs -0.2 cm {\bf .} ---  }
{  \end{remarq}}

\newtheorem{preuve}{\sc{Preuve}}

\def \vs {\vskip}
\def \hs {\hskip}
\def \noi {\noindent}
\def \fin {\hfill$\Box$}

\def \cad {c'est-{\`a}-dire }

\def \oo {{\cal O}}

\def \a {{\alpha}}

\def \ga {{\gamma}}

\def \vp {\varphi}

\def \p {{\mathbb P}}
\def \Z {{\mathbb{Z}}}

\def \G {{\mathbb G}}

\def \fll {\longrightarrow}

\def \Hom {{\rm Hom}}

\def \pic {{\rm Pic}}
\def \Stab {{\rm Stab}}

\def \pu {{\mathbb P}^1}

\def \MorC #1#2{{\bf{Hom}}_{#1}(C,#2)}

\def \sca #1#2{\left\langle#1,#2\right\rangle}

\def \Xt {{\widetilde{X}}}
\def \Vt {{\widetilde{V}}}
\def \Ut {{\widetilde{U}}}
\def \ft {{\widetilde{f}}}

\def \at {{\widetilde{\a}}}

\def \xih {{\widehat{\xi}}}

\def \cC {{\mathfrak{C}}}

\def \cw {\overline{w}}
\def \cs {\overline{s}}

\def \app {A_1^+(\Xt)}

\def \wp {{W_{\bullet}}}

\begin{document}

~
\vs -1 cm

\centerline{\large{\uppercase{\bf{Courbes elliptiques sur la vari{\'e}t{\'e} }}}}
\centerline{\large{\uppercase{\bf{spinorielle de
        dimension 10}}}}

\vs 0.5 cm

\centerline{\Large{Nicolas \textsc{Perrin}}}

\vs 1 cm

\centerline{\Large{\bf Introduction}} 

\vs 0.3 cm

Dans cet article, nous r{\'e}pondons {\`a} une question de Dimitri  Markushevich
concernant l'ir\-r{\'e}ducti\-bilit{\'e} du sch{\'e}ma des courbes elliptiques sur la
vari{\'e}t{\'e} spinorielle de dimension 10.

Soit $G$ le groupe $S0(10)$ et $P$ le parabolique maximal associ{\'e} {\`a}
la racine simple $\a_5$ avec les notation de \cite{bourb}. La vari{\'e}t{\'e}
homog{\`e}ne $G/P$ correspond {\`a} la grassmannienne des sous-espaces
totalement isotropes maximaux (ici de dimension 5) de type 1 dans un
espace vectoriel de dimension 10 muni d'une forme quadratique non
d{\'e}g{\'e}n{\'e}r{\'e}e. Elle est lisse de dimension 10, nous la noterons $X$. Nous
montrons le r{\'e}sultat suivant :

\begin{theo}
  Soit $C$ une courbe elliptique lisse et soit $\a\in A_1(X)$ une
  classe de 1-cycles sur $X$ de degr{\'e} $d$. Si $d\geq4$, alors le
  sch{\'e}ma $\MorC{\a}{X}$ des morphismes de $C$ vers $X$ de classe $\a$
  est irr{\'e}ductible de dimension $8d$.
\end{theo}

Pour montrer ce r{\'e}sultat, nous utilisons la r{\'e}solution de
Bott-Samelson $\pi:\Xt_\wp\to X$ associ{\'e}e {\`a} un drapeau complet $\wp$
(cf. paragraphe \ref{Bott}). Contrairement au cas des courbes
rationnelles que nous avons trait{\'e} dans \cite{perrin} ou
\cite{Perrin}, nous avons besoin d'imposer des conditions $(*)$
(cf. paragraphe \ref{irredu}) sur les classes $\at\in A_1(\Xt)$ telles
que $\pi_*(\at)=\a$ pour obtenir l'irr{\'e}ductibilit{\'e} du sch{\'e}ma
$\MorC{\at}{\Xt_\wp}$. Nous montrons (propositions \ref{intersection} et
\ref{reste}) que pour (essentiellement) tout morphisme
$f\in\MorC{\a}{X}$, il existe un drapeau complet $\wp$ tel que $f$ se
rel{\`e}ve dans $\Xt_\wp$ en un morphisme $\ft$ satisfaisant les
conditions $(*)$. On conclue alors {\`a} l'irr{\'e}ductibilit{\'e} en faisant
varier $\wp$.

\begin{rema}
(\i) Pour les morphismes de degr{\'e} 2 et 3 
on montre facilement que les vari{\'e}t{\'e}s sont  irr{\'e}ductibles de dimensions
respectives 19 et 25.

(\i\i) La technique pr{\'e}sent{\'e}e ici doit pouvoir se g{\'e}n{\'e}raliser aux
vari{\'e}t{\'e}s homog{\`e}nes minuscules. 
C'est ce que nous faisons dans \cite{PEcarquois} pour les groupes
  classiques. 
\end{rema}


\section{R{\'e}solution de Bott-Samelson}
\label{Bott}

Nous reprendrons la plupart des notations de \cite{Perrin}.
Soit $W$ le groupe de Weyl de $G$ et $W_P$ le sous-groupe de $W$
stabilisant $P$. On note $w_0$ l'{\'e}l{\'e}ment de longueur maximal dans $W$
et $\cw_0$ sa classe dans $W/W_P$. On a une {\'e}criture r{\'e}duite de
$w_0$ et m{\^e}me de $\cw_0$ :
$$\cw_0=\cs_4\ \!\cs_3\ \!\cs_2\ \!\cs_5\ \!\cs_1\ \!
\cs_3\ \!\cs_2\ \!\cs_4\ \!\cs_3\ \!\cs_5$$
\vs -0.9 cm \hfill$(*)$
\vs 0.2 cm
\noi
o{\`u} $\cs_i$ est la classe dans $W/W_P$ de la sym{\'e}trie par rapport {\`a} la
racine simple $\a_i$ (toujours avec les notations de \cite{bourb}).

Une fois cette d{\'e}composition fix{\'e}e, {\`a} tout drapeau complet $\wp$, on
peut associer (cf. \cite{Demazure}) une
vari{\'e}t{\'e} de Bott-Samelson $\Xt_\wp$. Cette vari{\'e}t{\'e} peut {\^e}tre
d{\'e}crite comme une suite de fibrations
$$\Xt_\wp=X_{10}\stackrel{f_{10}}{\fll}\cdots
\stackrel{f_{2}}{\fll}X_1\stackrel{f_{1}}{\fll}X_0\simeq{\rm
  Spec}(k).$$ 
o{\`u} les $f_i$ sont des fibrations en droite projectives. Chacune de ces
fibrations $f_i$ est munie d'une section $\sigma_i$ et on note $\xi_i$
le diviseur de $X_i$ donn{\'e} par $\sigma_i(X_{i-1})$. Par abus de
notations, nous noterons encore $\xi_i$ l'image r{\'e}ciproque de ce
diviseur dans $\Xt_\wp$. Nous noterons $T_i$ le fibr{\'e} tangent relatif
de la fibration $f_i$ et $p_i$ le morphisme de $\Xt_\wp$ vers $X_i$.

%
%
%
%


%
%
%
%
%

\section{Fibrations}

\subsection{Irr{\'e}ductibilit{\'e}}
\label{irredu}

Nous montrons une proposition qui permet de remonter
l'irr{\'e}ductibilit{\'e} du sch{\'e}ma des morphismes {\`a} travers les fibrations en
droites projectives. Cependant, contrairement au cas des courbes
rationnelles (cf. \cite{perrin} prop. 4), la seule condition d'avoir
un degr{\'e} relatif positif ne suffit plus. Soit $C$ une courbe
elliptique lisse.

\begin{prop}
\label{irred}
  Soit $\vp:X\to Y$ une fibration en droites projectives munie d'une
  section $\sigma$ et soit $\at\in A_1(X)$ une classe de
  1-cycles. Notons $T$ le fibr{\'e} tangent relatif et $\xi$ le diviseur
  donn{\'e} par la section. 
Supposons que $\at$ v{\'e}rifie $\at\cdot\xi\geq0$ et $\at\cdot(T-\xi)>0$.

Si $\MorC{\vp_*\at}{Y}$ est irr{\'e}ductible, alors $\MorC{\at}{X}$ l'est
et  
$$\dim(\MorC{\at}{X})=\dim(\MorC{\vp_*\at}{Y})+\at\cdot T.$$
\end{prop}

\dm
Notons $E$ un fibr{\'e} vectoriel de rang 2 sur $Y$ tel que
$X=\p_Y(E)$. La section $\sigma$ est donn{\'e}e par une surjection
$E\to L$ o{\`u} $L$ est inversible. Notons $N$ le fibr{\'e} inversible noyau
de cette surjection. 

Nous {\'e}tudions la fibre du morphisme $\MorC{\at}{X}\to\MorC{\vp_*\at}{Y}$
au dessus de la fl{\`e}che $f:C\to Y$. Un {\'e}l{\'e}ment de la fibre est donn{\'e}
par un rel{\`e}vement de $f$, c'est {\`a} dire par une surjection $f^*E\to M$
o{\`u} $M$ est inversible sur $C$ avec $2\deg(M)-\deg(f^*E)=\at\cdot
T$. Un {\'e}l{\'e}ment de la fibre est donc donn{\'e} par un fibr{\'e} inversible $M$
de degr{\'e} $d=\frac{\deg(f^*E)+\at\cdot T}{2}$ et par un {\'e}l{\'e}ment
surjectif de $\p\Hom(f^*E,M)$. On a $\at\cdot\xi=\deg(M)-\deg(f^*N)$
et $\at\cdot(T-\xi)=\deg(M)-\deg(f^*L)$ et on doit distinguer deux
cas.

Si $\at\cdot\xi>0$ alors $\Hom(f^*E,M)$ est la somme directe de
$\Hom(f^*N,M)$ et de $\Hom(f^*L,M)$ et est de dimension contante (par
rapport {\`a} $f$) {\'e}gale {\`a} $\at\cdot T$ et le choix de $M$ est
libre. La fibre est donn{\'e}e par le choix de $M$ puis d'une surjection
$f^*E\to M$ \cad par un ouvert (donn{\'e} par la condition de
surjectivit{\'e}) non vide de $\p\Hom(f^*E,M)\times\pic_{d}(C)$. On a donc
une fibration lisse de dimension $\at\cdot T$ au-dessus de
$\MorC{\vp_*\at}{Y}$ d'o{\`u} le r{\'e}sultat.

Si $\at\cdot\xi=0$ alors si $M\not\simeq f^*N$ on a $\Hom(f^*N,M)=0$
donc toute fl{\`e}che $f^*E\to M$ se factorise par $f^*L$ et on ne peut
avoir de fl{\`e}che surjective car
$\deg(M)-\deg(f^*L)=\at\cdot(T-\xi)>0$. Pour tout {\'e}l{\'e}ment de la fibre,
on doit donc avoir un isomorphisme $M\simeq f^*N$. Mais alors comme
$\deg(M)-\deg(f^*L)=\at\cdot(T-\xi)>0$, on a $E=M\oplus f^*L$ et
$\Hom(f^*E,M)$ est la somme directe de $\Hom(f^*N,M)$ et de
$\Hom(f^*L,M)$ et est de dimension contante (par rapport {\`a} $f$) {\'e}gale
{\`a} $\at\cdot T+1$. La fibre est donn{\'e}e par un ouvert non vide
de $\p\Hom(f^*E,M)$. 
On a donc une fibration lisse de dimension $\at\cdot T$ au-dessus de
$\MorC{\vp_*\at}{Y}$.\fin

\begin{rema}
%
On sait calculer les $T_i$ en fonction des $\xi_i$ (cf. \cite{Perrin}
prop. 2.11) :
$$T_i=\sum_{k=1}^i\sca{\ga_k^\vee}{\ga_i}\xi_i$$
o{\`u} les $\ga_i$ sont des racines positives d{\'e}finies {\`a} partir de
la d{\'e}composition  r{\'e}duite $(*)$ (cf.\cite{Perrin}).
On peut ici calculer tous les $\sca{\ga_k^\vee}{\ga_i}$ et on a (le
r{\'e}sultat est sym{\'e}trique en $k$ et $i$ car $R=R^\vee$) :
\begin{small}
$$
\begin{tabular}{|c||c|c|c|c|c|c|c|c|c|c|}
\hline
&1&2&3&4&5&6&7&8&9&10\\
\hline
\hline
1&2&1&1&1&1&1&1&0&0&0\\
\hline
2&&2&1&1&1&0&0&1&1&0\\
\hline
3&&&2&0&1&1&0&1&0&1\\
\hline
4&&&&2&0&1&1&1&1&0\\
\hline
5&&&&&2&0&1&0&1&1\\
\hline
6&&&&&&2&1&1&0&1\\
\hline
7&&&&&&&2&0&1&1\\
\hline
8&&&&&&&&2&1&1\\
\hline
9&&&&&&&&&2&1\\
\hline
10&&&&&&&&&&2\\
\hline
\end{tabular}$$

\end{small}
\vs -0.45 cm\hs 11.2 cm .

\vs 0.45 cm
On peut lire l'{\'e}criture de $T_i$ sur la $i^{i\grave{e}me}$
colonne. En particulier, on voit que si une classe $\at$ est
strictement positive par rapport {\`a} $\xi_1$, $\xi_2$ et $\xi_3$
et positive par rapport {\`a} tous les $\xi_i$, alors elle est strictement
positive par rapport {\`a} tous les fibr{\'e}s $T_i-\xi_i$. Notons $\app$
l'ensemble des classes $\at$ qui sont strictement positives par rapport
{\`a} $\xi_1$, $\xi_2$ et $\xi_3$ et positives par rapport {\`a} tous les $\xi_i$. 
\end{rema}

\subsection{Calcul de la dimension}
\label{dimension}

La proposition pr{\'e}c{\'e}dente nous permet donc de montrer que pour tout
$\a\in\app$, le sch{\'e}ma $\MorC{\a}{\Xt}$ est irr{\'e}ductible de dimension
$$\at\cdot\sum_{i=1}^{10}T_i=
\at\cdot(8\xi_1+7\xi_2+6\xi_3+6\xi_4+5\xi_5+5\xi_6+4
\xi_7+4\xi_8+3\xi_9 +2\xi_{10}).$$
Le g{\'e}n{\'e}rateur ample du groupe de Picard de $X$ se rel{\`e}ve dans $\Xt$ en
un diviseur $\xih$ qui s'{\'e}crit
$$\xih=\sum_{i=1}^{10}\xi_i.$$
Ainsi si $\at$ est une classe de $\app$
telle que $\deg(\pi_*(\at))=\at\cdot\xih=d$, alors
le sch{\'e}ma $\MorC{\at}{\Xt}$ est irr{\'e}ductible de dimension :
$$8d-\at\cdot\xi_2-2\at\cdot\xi_3-2\at\cdot\xi_4-3
\at\cdot\xi_5-3\at\cdot\xi_6-4
\at\cdot\xi_7-4\at\cdot\xi_8-5\at\cdot\xi_9 -6\at\cdot\xi_{10}.$$
En particulier la dimension est maximale exactement lorsque
$\at\cdot\xi_i=0$ pour $i$ diff{\'e}rent de 1, 2 et 3 et $\at\cdot\xi_i=1$ pour
$i\in\{1;2;3\}$. La dimension est alors $8d-3$.

\section{Description de $\Xt_\wp$}
\label{diviseurs}

\subsection{Configurations}

Nous d{\'e}crivons ici la r{\'e}solution de Bott-Samelson sous la forme d'une
vari{\'e}t{\'e} de configuration. Ceci a {\'e}t{\'e} r{\'e}alis{\'e} par P. Magyar
\cite{Magyar}. Nous avons ici
$$\Xt_\wp=\left\{(\Vt_1,\Vt_2,\Vt_2',\Vt_3,\Vt_3',\Vt_3'',\Vt_4,
  \Vt_4',\Vt_5,\Vt_5')\ 
/\ (**)\right\}$$
o{\`u} l'indice $i$ pour $i\in\{1;2;3\}$ indique qu'on a un sous-espace
  totalement isotrope de dimension $i$ et o{\`u} l'indice 4 respectivement
  5 indique qu'on a un sous-espace totalement isotrope maximal de type
  2 ou 1 respectivement. Les conditions $(**)$ sont donn{\'e}es par les 10
  conditions :
$$\begin{array}{ccccc}
\Vt_4\supset W_3,&W_2\subset\Vt_3\subset
W_5\cap\Vt_4,&W_1\subset\Vt_2\subset
\Vt_3,&\Vt_3\subset\Vt_5,&\Vt_1\subset \Vt_2,\\
\Vt_2\subset\Vt_3'\subset \Vt_5\cap\Vt_4,&\Vt_1\subset\Vt_2'\subset
\Vt_3',&\Vt_3'\subset\Vt_4',&\Vt_2'\subset\Vt_3''\subset
\Vt_5\cap\Vt_4',&\Vt_3''\subset\Vt_5'.\\
\end{array}$$

On peut en fait "visualiser" cette configuration sous la forme d'un
carquois (voir \cite{PEcarquois} pour plus de d{\'e}tails) o{\`u} une fl{\`e}che
signifie qu'on a une inclusion dans le seul sens possible. On ne fait
pas appara{\^\i}tre les points du drapeau de d{\'e}part $\wp$. Ainsi dans le
cas pr{\'e}sent le carquois est 
$$
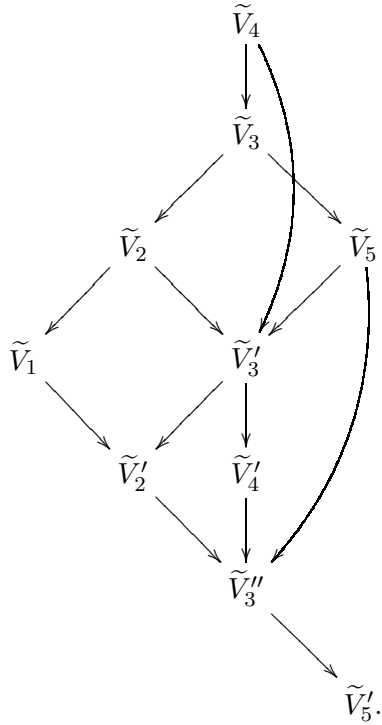
\xymatrix{
&&\Vt_4\ar[d]\ar@/^1.5pc/[ddd]&\\
&&\Vt_3\ar[dl]\ar[dr]&\\
&\Vt_2\ar[dl]\ar[dr]&&\Vt_5\ar[dl]\ar@/^1.5pc/[dddl]\\
\Vt_1\ar[dr]&&\Vt_3'\ar[d]\ar[dl]&\\
&\Vt_2'\ar[dr]&\Vt_4'\ar[d]&\\
&&\Vt_3''\ar[dr]&\\
&&&\Vt_5'.}$$
Cette {\'e}criture sous forme de carquois n'appara{\^\i}t pas explicitement
dans \cite{Magyar}. Les liens entre d{\'e}composition r{\'e}duite et carquois
plong{\'e} dans le complexe de Coxeter sont mis en {\'e}vidence par
S. Zelikson \cite{Zelikson} ou dans le mod{\`e}le des galeries
combinatoires de C.E. Contou-Carr{\`e}re \cite{CC}, voir aussi S. Gaussent
\cite{Gaussent}. 

Les fibrations $f_i$ sont donn{\'e}es par les morphismes d'oubli en
remontant les fl{\`e}ches du carquois : on oublie $\Vt_5'$ puis $\Vt_3''$
et ainsi de suite. Les sections sont alors donn{\'e}es par les {\'e}galit{\'e}s 
$$\begin{array}{|c|c|c|c|c|}
\hline
\xi_1:\Vt_4=W_4,&\xi_2:\Vt_3=W_3,&\xi_3:\Vt_2=W_2,&\xi_4:\Vt_5=W_5,
&\xi_5:\Vt_1=W_1,\\ 
\hline
\xi_6:\Vt_3'=\Vt_3,&\xi_7:\Vt_2'=\Vt_2,&
\xi_8:\Vt_4'=\Vt_4,&\xi_9:\Vt_3''=\Vt_3',&\xi_{10}:\Vt_5'=\Vt_5.\\ 
\hline
\end{array}$$

Le morphisme de $\Xt_\wp$ vers $X$ est donn{\'e} par la derni{\`e}re
projection \cad par l'espace $\Vt'_5$. Si on a un point $V\in X$ en
position g{\'e}n{\'e}rale par rapport au drapeau $\wp$, alors on peut
retrouver les espaces de la configuration pr{\'e}c{\'e}dente gr{\^a}ce aux
formules :
$$\begin{array}{|c|c|c|c|c|}
\hline
\Vt'_5=V&\Vt_1=V\cap
W_5&\Vt_2=\Vt_1+W_1&\Vt_3=\Vt_1+W_2&\Vt_4\supset\Vt_1+W_3\\
\hline
\Vt_2'=V\cap
\Vt_4&\Vt_3'=\Vt_2'+\Vt_2&\Vt_5\supset\Vt_2'+\Vt_3&\Vt_3''=
V\cap\Vt_5&\Vt_4'\supset\Vt_3'+\Vt_3''\\
\hline
\end{array}$$
en remarquant que lorsqu'on a un sous-espace totalement isotrope
maximal de type fix{\'e}, il est d{\'e}termin{\'e} par un sous-espace totalement
isotrope de dimension 4 qu'il contient.

En particulier, si un morphisme $f:C\to X$ rencontre l'ouvert des
points en position g{\'e}n{\'e}rale par rapport {\`a} $\wp$, on peut d{\'e}finir sur
un ouvert de $C$ et donc par prolongement sur $C$ tout enti{\`e}re une
section $\ft$ de $f$ dans $\Xt_\wp$ donn{\'e}e par les formules du tableau
pr{\'e}c{\'e}dent.

\subsection{Images par $\pi$}

Il est facile de d{\'e}crire l'image des diviseurs $\xi_i$ dans $X$. On a 

\vs 0.2 cm


\centerline{\begin{tabular}{cc}
$\pi(\xi_1)=\{V\in X\ /\ \dim(V\cap W_4)\geq1\}$&$\pi(\xi_2)=\{V\in X\
/\ \dim(V\cap W_3)\geq1\}$\\
&\\
$\pi(\xi_3)=\{V\in X\ /\ \dim(V\cap W_2)\geq1\}$&
$\pi(\xi_4)=\{V\in X\ /\ \dim(V\cap W_5)\geq3\}$\\
&\\
$\pi(\xi_5)=\{V\in X\ /\ \dim(V\cap W_1)\geq1\}$&$\pi(\xi_6)=\left\{V\in X\ /\ 
\begin{array}{c}
\dim(V\cap W_2)\geq1\\
\dim(V\cap
W_5)\geq3
\end{array}\right\}$\\
&\\
$\pi(\xi_7)=\left\{V\in X\ /\ \begin{array}{c}
\dim(V\cap W_1)\geq1\\
\dim(V\cap
W_5)\geq3
\end{array}\right\}$&
$\pi(\xi_8)=\{V\in X\ /\ \dim(V\cap W_3)\geq2\}$\\
&\\
$\pi(\xi_9)=\left\{V\in X\ /\ \begin{array}{c}
\dim(V\cap W_1)\geq1\\
\dim(V\cap
W_3)\geq2
\end{array}\right\}$&$\pi(\xi_{10})=\{V\in X\ /\ \dim(V\cap
W_2)\geq2\}$.\\
\end{tabular}}


\section{{\'E}tude d'un ouvert de $X$}

Fixons un sous-espace totalement isotrope $W_5$ de dimension 5 et de
type 1. Si $\Stab(W_5)$ est le stabilisateur de $W_5$ dans $G$, les
orbites de ce stabilisateur sont donn{\'e}es par
$$\{V\in G/P\ /\ \dim(V\cap W_5)=1\},\ \ \{V\in G/P\ /\ \dim(V\cap
W_5)=3\}\ \ {\rm et}\ \ \{W_5\}.$$
La premi{\`e}re orbite est dense la seconde de codimension 3 et la dern{\`e}re de
codimension 10. Notons $U$ la premi{\`e}re orbite qui est un ouvert dense
de $X$ et $i$ l'inclusion. Le fait que $X\setminus U$ est de
codimension 3 et la proposition 2 de \cite{perrin} (qui
s'adapte 
au cas d'une courbe quelconque) montre que pour
{\'e}tudier l'irr{\'e}ductibilit{\'e} de $\MorC{\a}{X}$, il suffit d'{\'e}tudier celle
de $\MorC{i^*\a}{U}$. 

Sur l'ouvert $U$ on a un morphisme naturel 
$$\xymatrix{p_{W_5}:U\ar[r]&\p(W_5^\vee)\\
\ \ \ \ \ \ \ V\ar@{|->}[r]& [V\cap W_5].}$$
La proposition 5 de \cite{perrin} d{\'e}crit ce morphisme comme un fibr{\'e}
vectoriel au dessus de $\p(W_5^\vee)$ et la description donn{\'e}e dans
cette proposition permet d'identifier ce fibr{\'e} qui est
$\Lambda^2T_{\p(W_5^\vee)}(-1)$. Nous montrons maintenant la

\begin{prop}
  \label{intersection}
Soit $f:C\to U$ un morphisme, si $p_{W_5}\circ f(C)_{red}$ l'image
r{\'e}duite de $C$ par la compos{\'e}e $p_{W_5}\circ f$ n'est pas contenue
dans un plan, alors on peut compl{\'e}ter $W_5$ en un drapeau complet
$\wp$ tel que $f$ se rel{\`e}ve en $\ft$ dans $\Xt_\wp$ avec
$\ft(C)\cdot\xi_i\geq0$ pour tout $i$ et
$$\ft(C)\cdot\xi_1>0,\ \ft(C)\cdot\xi_2>0\ {\rm et}\
\ft(C)\cdot\xi_3>0.$$
\end{prop}

\dm
%
Soit $x_1\in C$ tel que $p_{W_5}(x_1)$ est un point lisse de
$p_{W_5}(C)_{red}$. Comme cette courbe n'est pas plane (et a fortiori
pas contenue dans une droite), il existe un point $x_2\in C$ tel que
$p_{W_5}(x_2)$ est un point lisse et n'est pas contenu dans la
tangente $\Theta_1$ {\`a} $p_{W_5}(C)_{red}$ en $p_{W_5}(x_1)$. Par
ailleurs, 
on peut quitte {\`a} changer $x_1$ choisir $x_2$ tel que la tangente
$\Theta_2$ {\`a} $p_{W_5}(C)_{red}$ en $p_{W_5}(x_2)$ ne passe pas par
$p_{W_5}(x_1)$. Comme la courbe $p_{W_5}(C)_{red}$ n'est pas plane, il
existe un point $x_3\in C$ tel que $p_{W_5}(x_3)$ n'est contenu dans
aucun des deux plans $\Pi_1$ et $\Pi_2$ respectivement engendr{\'e}s par
la droite $\Theta_i$ et le point $p_{W_5}(x_{3-i})$ pour
$i\in\{1;2\}$. De m{\^e}me, il existe un point $x_4\in C$ tel que
$p_{W_5}(x_4)$ n'est pas contenu dans le plan $\Pi$ engendr{\'e} par
$p_{W_5}(x_1)$, $p_{W_5}(x_2)$ et $p_{W_5}(x_3)$ ni dans les plans
$\Pi_1$ et $\Pi_2$.

L'intersection du plan $\Pi$ et des plans $\Pi_1$ et $\Pi_2$ est la
droite $(p_{W_5}(x_1),p_{W_5}(x_2))$ qui rencontre les droites
$\Theta_i$ aux points $p_{W_5}(x_i)$ pour $i\in\{1;2\}$. Ainsi le plan
$\Pi$ ne contient pas $\Theta_i$ pour $i\in\{1;2\}$ et nous pouvons
choisir un sous-espace vectoriel $U_4$ de dimension 4 de $W_5$
contenant $p_{W_5}(x_1)$, $p_{W_5}(x_2)$ et $p_{W_5}(x_3)$ mais ne
contenant pas $p_{W_5}(x_4)$ ni $\Theta_i$  pour $i\in\{1;2\}$. Nous
choisissons ensuite un sous-espace vectoriel $W_3$ de dimension 3 de
$U_4$ contenant $p_{W_5}(x_1)$ et $p_{W_5}(x_2)$ mais ne contenant pas
$p_{W_5}(x_3)$, un sous-espace vectoriel $W_2$ de dimension 2 de $W_3$
contenant $p_{W_5}(x_1)$ mais ne contenant pas $p_{W_5}(x_2)$ et enfin
un sous-espace vectoriel $W_1$ de dimension 1 de $W_2$ ne contenant
pas $p_{W_5}(x_1)$. Pour compl{\'e}ter le drapeau, nous prenons pour $W_4$
le sous-espace totalement isotrope de dimension 5 de type 2 contenant
$U_4$. Nous obtenons ainsi un drapeau $\wp$. Il nous suffit pour
terminer de monter le 

\begin{lemm}
  Le morphisme $f:C\to X$ se rel{\`e}ve en $\ft$ dans $\Xt_\wp$ avec
  $\ft(C)\cdot\xi_i\geq0$ pour tout $i$ et 
$$\ft(C)\cdot\xi_1>0,\ \ft(C)\cdot\xi_2>0\ {\rm et}\
\ft(C)\cdot\xi_3>0.$$
\end{lemm}

\dm
Pour montrer que $f$ se rel{\`e}ve en $\ft$ dans $\Xt_\wp$ avec
$\ft(C)\cdot\xi_i\geq0$ pour tout $i$, il suffit de montrer qu'un
point de $f(C)$ est en position g{\'e}n{\'e}rale par rapport au drapeau
$\wp$. C'est le cas du point $f(x_4)$. En effet, on sait que
l'intersection de $f(x_4)$ avec $W_5$ est $p_{W_5}(x_4)$ donc $f(x_4)$
ne peut rencontrer $W_1$, $W_2$, $W_3$ ni $U_4$. Par ailleurs, si
$f(x_4)$ rencontrait $W_4$, cette intersection serait de dimension au
moins 2 et rencontrerait donc $U_4$, c'est impossible.

Il reste {\`a} montrer que $\ft(C)$ rencontre $\xi_i$ en $f(x_{4-i})$ pour
$i\in\{1;2;3\}$. Le m{\^e}me raisonnement que ci-dessus implique que
$f(x_3)$ ne rencontre pas $W_1$, $W_2$ ni $W_3$ mais rencontre $W_5$
en dimension 1 exactement et rencontre {\'e}galement $W_4$. Ainsi $f(x_3)$
est contenu dans $\pi(\xi_1)$ mais dans aucun des $\pi(\xi_i)$ pour
$i>1$. Le point $\ft(x_3)$ est donc contenu dans $\xi_1$.

De la m{\^e}me mani{\`e}re, on voit que $f(x_2)$ ne rencontre pas $W_1$ ni
$W_2$  mais rencontre $W_3$ en dimension 1 exactement. Ainsi $f(x_2)$
est contenu dans $\pi(\xi_1)$ et $\pi(\xi_2)$ mais dans aucun des
$\pi(\xi_i)$ pour $i>2$. Par ailleurs on a une application $\Ut_4$
d{\'e}finie sur un ouvert de $C$ par 
$$\Ut_4(x)=p_{W_5}(x)+W_3$$ 
qui se prolonge sur $C$ tout enti{\`e}re. L'image de $x_2$ est donn{\'e}e par
$\Ut_4(x_2)=\Theta_2+W_3\neq U_4$. Ainsi l'espace $\Vt_4(x_2)$ qui
permet de d{\'e}finir $\ft(x_2)$ est diff{\'e}rent de $W_4$ et
$\ft(x_2)\not\in\xi_1$. Le point $\ft(x_2)$ est donc contenu dans
$\xi_2$.

Enfin, on voit que $f(x_1)$ ne rencontre pas $W_1$ mais rencontre
$W_2$ en dimension 1 exactement. Ainsi $f(x_1)$ est contenu dans
$\pi(\xi_1)$ $\pi(\xi_2)$ et $\pi(\xi_3)$ mais dans aucun des
$\pi(\xi_i)$ pour $i>3$. Par ailleurs on a deux applications $\Vt_3$
et $\Ut_4$ d{\'e}finies sur un ouvert de $C$ par  
$$\Vt_3(x)=p_{W_5}(x)+W_2\ \ \ {\rm et}\ \ \ \Ut_4(x)=p_{W_5}(x)+W_3$$
qui se prolongent sur $C$ tout enti{\`e}re. Les images de $x_1$ sont
donn{\'e}es par $\Vt_3(x_3)=\Theta_1+W_2\neq W_3$ et
$\Ut_4(x_1)=\Theta_1+W_3\neq U_4$. Ainsi on a $\ft(x_1)\not\in\xi_1$
et $\ft(x_1)\not\in\xi_2$. Le point $\ft(x_1)$ est donc contenu dans
$\xi_3$.\fin

\begin{rema}
  \label{fibrep}
  La construction pr{\'e}c{\'e}dente montre que pour un morphisme $f:C\to U$
  donn{\'e} tel que $p_{W_5}\circ f(C)_{red}$ n'est pas contenue dans un
  plan, les drapeaux qui compl{\`e}tent $W_5$ en un drapeau complet
  v{\'e}rifiant les conditions de la  proposition \ref{intersection}
  forment une vari{\'e}t{\'e} de dimension 7 : choix des points $x_i$ pour
  $i\in\{1;2;3\}$ (dimension 3) puis choix des sous-espaces $W_i$ pour 
  $i\in\{1;2;3;4\}$ (dimension 1 pour chaque $W_i$).
\end{rema}


Nous montrons maintenant que les morphismes qui ne satisfont pas les
conditions de la proposition \ref{intersection} ne peuvent former une
composante irr{\'e}ductible du sch{\'e}ma des morphismes. Rappelons que le
morphisme $p_{W_5}:U\to\p(W_5^\vee)$ est le morphisme de projection du
fibr{\'e} vectoriel $\Lambda^2T_{\p(W_5^\vee)}(-1)$.

\begin{prop}
\label{reste}
  Soit $\a\in A_1(X)$ une classe de degr{\'e} $d\geq2$. Le ferm{\'e} du sch{\'e}ma des
  morphismes $\MorC{i^*\a}{U}$ form{\'e} des fl{\`e}ches $f$ telles que la
  courbe $p_{W_5}\circ f(C)_{red}$ est contenue dans un plan est de
  dimension au plus $6d+7$. 

Lorsque $d\geq4$, ce ferm{\'e} ne peut former une
composante irr{\'e}ductible de $\MorC{i^*\a}{U}$.
\end{prop}

\dm
Nous consid{\'e}rons le morphisme
$\MorC{i^*\a}{U}\to\MorC{{p_{W_5}}_*i^*\a}{\p(W_5^\vee)}$ induit par
$p_{W_5}$. Pour d{\'e}finir ${p_{W_5}}_*$ nous consid{\'e}rons les classes de
1-cycles comme des {\'e}l{\'e}ments du dual du groupe de Picard et utilisons
la transpos{\'e}e de ${p_{W_5}}^*$. L'image du ferm{\'e} consid{\'e}r{\'e} dans
$\MorC{{p_{W_5}}_*i^*\a}{\p(W_5^\vee)}$ est irr{\'e}ductible de dimension
$3d+6$ et contient le ferm{\'e} des courbes dont le support r{\'e}duit est
une droite qui est irr{\'e}ductible de dimension $2d+6$.

La fibre du morphisme
$\MorC{i^*\a}{U}\to\MorC{{p_{W_5}}_*i^*\a}{\p(W_5^\vee)}$ au dessus de
$f:C\to\p(W_5^\vee)$ est donn{\'e}e par
$H^0(f^*\Lambda^2T_{\p(W_5^\vee)}(-1))$. Pour un morphisme dont
l'image est contenue dans un plan, on a 
$$f^*T_{\p(W_5^\vee)}(-1)=\oo_C^2\oplus E$$
o{\`u} $E$ est un fibr{\'e} de rand 2 et de degr{\'e} $d$ sur $C$ engendr{\'e} par ses
sections (car $T_{\p(W_5^\vee)}(-1)$ l'est). De plus $E$ a un facteur
trivial si et seulement si l'image de $f$ est contenue dans une
droite. Ainsi on a 
$$f^*\Lambda^2T_{\p(W_5^\vee)}(-1)=\oo_C\oplus\oo_C^2\otimes
E\oplus\Lambda^2E.$$ 
Si $E$ n'a pas de facteur trivial \cad lorsque l'image de $f$ n'est
pas contenue dans une droite, alors  la dimension du groupe
$H^0(f^*\Lambda^2T_{\p(W_5^\vee)}(-1))$ est $3d+1$. Si $E$ a un
facteur trivial (n{\'e}cessairement unique) \cad lorsque l'image de $f$
est contenue dans une droite, alors la dimension du groupe
$H^0(f^*\Lambda^2T_{\p(W_5^\vee)}(-1))$ est $3d+3$. Ainsi le ferm{\'e}
recherch{\'e} a une dimension au plus {\'e}gale {\`a} $\max(6d+7;5d+9)=6d+7$ pour
$d\geq2$.

Enfin la dimension d'une composante irr{\'e}ductible de $\MorC{i^*\a}{U}$
est au moins $8d$ qui est strictement sup{\'e}rieur {\`a} $6d+7$ d{\`e}s que
$d\geq4$ ce qui prouve la derni{\`e}re assertion.\fin

\section{Irr{\'e}ductibilit{\'e}}

\subsection{L'incidence}

Soit $B$ un sous-groupe de Borel de $G$, la vari{\'e}t{\'e} des drapeaux
complets est $G/B$. Les drapeaux completant $W_5$ en un drapeau
complet forment la vari{\'e}t{\'e} $\Stab(W_5)/B$ qui est de dimension
10. Nous avons vu qu'{\`a} chaque drapeau complet $\wp$ correspond une
r{\'e}solution de Bott-Samelson $\Xt_\wp\to X$. Consid{\'e}rons l'incidence
suivante :
$$
\xymatrix{
I  \ar@{->}[r]^-{p} \ar[d]_q & {\MorC{i^*\a}{U}}\\
\Stab(W_5)/B &
\\}
$$
o{\`u} $I$ est form{\'e}e des couples $(\wp,f)\in \Stab(W_5)/B\times\MorC{i^*\a}{U}$
v{\'e}rifiant les conclusions de la proposition \ref{intersection}. Nous
avons vu {\`a} la proposition \ref{reste} que la fl{\`e}che
$p:I\to\MorC{i^*\a}{U}$ est dominante d{\`e}s que $d\geq4$. Rappelons que
nous avons vu {\`a} la remarque \ref{fibrep} que la fibre de $p$ est de
dimension 7.

\subsection{La fibre de $q$}

Fixons $\wp$ un drapeau complet. La fibre de $q$ au dessus de $\wp$
est form{\'e} des morphismes $f\in\MorC{i^*\a}{U}$ qui se rel{\`e}vent en $\ft$
dans $\Xt_\wp$ et tels que $\ft(C)$ rencontre $\xi_1$, $\xi_2$ et
$\xi_3$.
On a donc une fl{\`e}che surjective
$$\coprod_{\at\in\app,\ \pi_*(\at)=\a}\MorC{\at}{\Xt_\wp}\to q^{-1}(\wp).$$
Nous pouvons faire cette construction en famille au dessus de
$\Stab(W_5)/B$.

Pour tout {\'e}l{\'e}ment $\at\in\app$ tel que la dimension de
$\MorC{\at}{\Xt_\wp}$ est strictement inf{\'e}rieure {\`a} $8d-3$ (\cad,
d'apr{\`e}s le paragraphe \ref{dimension} lorsque $\at$ ne v{\'e}rifie pas
$\at\cdot\xi_i=0$ pour $i$ diff{\'e}rent de 1, 2 et 3 et $\at\cdot\xi_i=1$
pour $i\in\{1;2;3\}$) la famille de morphismes $\MorC{\at}{\Xt_\wp}$
lorsque $\wp$ varie est de dimension strictement inf{\'e}rieure {\`a} $8d+7$
et son image dans $\MorC{i^*\a}{U}$ est strictement inf{\'e}rieure {\`a}
$8d$. Elle ne peut former une composante irr{\'e}ductible de
$\MorC{i^*\a}{U}$.

Ainsi si $\at_0$ est l'{\'e}l{\'e}ment de $\app$ tel que $\pi_*(\at)=\a$ et
$\at\cdot\xi_i=0$ pour $i$ diff{\'e}rent de 1, 2 et 3 et $\at\cdot\xi_i=1$
pour $i\in\{1;2;3\}$, la famille de morphismes $\MorC{\at_0}{\Xt_\wp}$
lorsque $\wp$ varie s'envoie de mani{\`e}re dominante sur
$\MorC{i^*\a}{U}$. Comme cette famille est irr{\'e}ductible gr{\^a}ce {\`a} la
proposition \ref{irred} on a l'irr{\'e}ductibilit{\'e} de $\MorC{\i^*\a}{U}$
puis gr{\^a}ce {\`a} la proposition 2 de \cite{perrin} : 

\begin{theo}
  Pour $\a\in A_1(X)$ tel que $d=\deg(\a)\geq4$, le sch{\'e}ma
  $\MorC{\a}{X}$ est irr{\'e}ductible de dimension $8d$.
\end{theo}


\begin{small}

\vs 0.2 cm

\noi
{\textsc{Institut de Math{\'e}matiques de Jussieu}}

\vs -0.1 cm

\noi
{\textsc{175 rue du Chevaleret}}

\vs -0.1 cm

\noi
{\textsc{75013 Paris,}} \hs 0.2 cm{\textsc{France.}}

\vs -0.1 cm

\noi
{email : \texttt{nperrin@math.jussieu.fr}}

\end{small}


\begin{thebibliography}{FMcPSS}
\bibitem[Bo]{bourb} \textit{Nicolas Bourbaki}: {\'E}l{\'e}ments de
    math{\'e}matique. Fasc. XXXIV. Groupes et alg{\`e}bres de Lie. Chapitre
    IV : Groupes de Coxeter et syst{\`e}mes de Tits.
   Chapitre V: Groupes engendr{\'e}s par des r{\'e}flexions. Chapitre VI:
   syst{\`e}mes de racines. Actualit{\'e}s Scientifiques et
   Industrielles, No. 1337 Hermann, Paris 1968.
\bibitem[CC]{CC} \textit{Carlos Enrique Contou-Carrere}: G{\'e}om{\'e}trie des
  groupes semi-simples, r{\'e}solutions {\'e}quivariantes et lieu singulier de
  leurs vari{\'e}t{\'e}s de Schubert. Th{\`e}se de l'universit{\'e} des sciences et
  techniques du languedoc (1983).
\bibitem[De]{Demazure} \textit{Michel Demazure}: D{\'e}singularisation des
    vari{\'e}t{\'e}s de Schubert g{\'e}n{\'e}ralis{\'e}es. Collection of articles
    dedicated to Henri Cartan on the occasion of his 70th birthday, I.
    Ann. Sci. {\'E}cole Norm. Sup. (4) 7 (1974), 53--88.  
\bibitem[Ga]{Gaussent} \textit{St{\'e}phane Gaussent}: The fibre of the
  Bott-Samelson resolution. Indag. Math. (N.S.) 12 (2001), no. 4,
  453--468.
\bibitem[Ma]{Magyar} \textit{Peter Magyar} : Schubert polynomials and
  Bott-Samelson varieties. Comment. Math. Helv. 73 (1998), no. 4,
  603--636.
\bibitem[P1]{perrin} \textit{Nicolas Perrin}: Courbes rationnelles sur
  les vari{\'e}t{\'e}s homog{\`e}nes. Annales de l'Institut  Fourier, 52, no.1
  (2002), pp 105-132.
\bibitem[P2]{Perrin} \textit{Nicolas Perrin}: Rational curves on
  minuscule Schubert varieties, pr{\'e}publication math.AG/0407123. 
\bibitem[P3]{PEcarquois} \textit{Nicolas Perrin}: Carquois et vari{\'e}t{\'e}s
  de Bott-Samelson, en pr{\'e}paration.
\bibitem[Ze]{Zelikson} \textit{Shmuel Zelikson} : Auslander-Reiten
  quivers and the Coxeter complex. Pr{\'e}publication math.QA/0208098 (2002).
  \end{thebibliography}
\end{document}